\documentclass[12pt,a4paper]{amsart}
\usepackage{amsfonts}
\usepackage{epsfig}
\usepackage{graphicx}
\usepackage{amsmath}
\usepackage{amssymb}
\usepackage{color}
\usepackage{mathrsfs}
\DeclareMathAlphabet{\mathpzc}{OT1}{pzc}{m}{it}

\newcommand{\Z}{\mathbb{Z}}

\DeclareMathAlphabet{\mathpzc}{OT1}{pzc}{m}{it}

\newcounter{main}

\newtheorem{theorem}{Theorem}[section]
\newtheorem{proposition}[theorem]{Proposition}
\newtheorem{lemma}[theorem]{Lemma}
\newtheorem{corollary}[theorem]{Corollary}

\newtheorem{remark}{Remark}[section]

\newtheorem{maintheorem}{Theorem}

\newcommand{\blanksquare}{\,\,\,$\sqcup\!\!\!\!\sqcap$}

\newcounter{example}
{{\stepcounter{example}}{\flushleft {\bf Example \arabic{example}:}}}%
{\par}

\title[Stable weakly shadowable volume-preserving systems]{Stable weakly shadowable volume-preserving systems are volume-hyperbolic}

\author[M. Bessa]{M\'{a}rio Bessa}

\address{Universidade da Beira Interior, Rua Marqu\^es d'\'Avila e Bolama,
  6201-001 Covilh\~a
Portugal.}
\email{bessa@fc.up.pt}

\author[M. Lee]{Manseob Lee}

\address{Department of Mathematics
Mokwon University, Daejeon 302-729, Republic of Korea.}
\email{lmsds@mokwon.ac.kr}

\author[S. Vaz]{Sandra Vaz}

\address{Universidade da Beira Interior, Rua Marqu\^es d'\'Avila e Bolama,
  6201-001 Covilh\~a
Portugal.}
\email{svaz@ubi.pt}

\date{\today}

\begin{document}

\begin{abstract}
We prove that any $C^1$-stably weakly shadowable vo\-lu\-me-preserving diffeomorphism defined on a compact manifold displays a dominated splitting $E\oplus F$. Moreover, both $E$ and $F$ are volume-hyperbolic. Finally, we prove the version of this result for divergence-free vector fields. As a consequence, in low dimensions, we obtain global hyperbolicity.
\end{abstract}

\maketitle

\tableofcontents

\section{Introduction}

\begin{subsection}{Shadowing in dynamical systems}
It is a very rich field in smooth dynamics the relation between the stability of a certain pro\-per\-ty (with respect to a given topology) and some hyperbolic behavior on the tangent action of the dynamical system. Structural stability, shadowing-type properties, robust transitivity, stable ergodicity, topological stability, expansiveness, specification, are some successful examples of that. Here we are interested in the weak shadowing property.

The shadowing in dynamics aims, in brief terms, to obtain shadowing of approximate trajectories in a given dynamical system by true orbits of the system. We refer the reader to the fairly complete exposition by Pilyugin's on the subject ~\cite{P}.

The weak shadowing property is a relaxed form of shadowing and first appear in the paper by Corless and Pilyugin (see~\cite{CP}) when related to the $C^0$-genericity of shadowing among dynamical systems. Informally speaking weak shadowing allows the pseudo-orbits to be approximated by true orbits if one forgets the time parameterization and consider only the distance between the orbit and the pseudo-orbit as two sets in the phase space. We intend, in this paper, to study the weak shadowing property for volume-preserving diffeomorphisms and also for volume-preserving flows.

There are limitations about the information we can capture from a fixed dynamical system that displays some shadowing-type property, since another system arbitrarily close to it may be absent of that property. Thence, it is of great utility and natural to consider that a selected model can be slightly perturbed in order to obtain the same property - the stable weakly shadowable systems. However, it is worth to mention that stability in the volume-preserving setting only allows us to consider perturbations which preserves the volume-form  and not evolving in the broader space of dissipative diffeomorphisms/flows. So the results already proved for dissipative diffeomorphisms/flows are not applicable to our conservative context.

In ~\cite{S2} it is proved that if a diffeomorphism defined in a surface has the $C^1$-stable weak shadowing property, then it satisfies the axiom A and the no-cycle condition. However, the converse does not hold (see \cite{P}). In \cite{PST} we can find more details on the relation between $C^1$-stability of weakly shadowing systems and structural stability in surfaces. In \cite{C},  the weak shadowing property is proved to be generic (in the $C^1$-sense) for diffeomorphisms in closed manifolds (this results is also valid in the volume-preserving context cf. ~\cite[\S 2.5]{C}).  In the two-dimensional volume-preserving case (thus symplectic because of the low dimension assumption), $C^1$-weakly shadowing implies hyperbolicity, \cite{L}, and  $C^1$-weakly shadowable symplectomorphisms are partially hyperbolic, \cite{BV}. In this paper we generalize the results in \cite{GSW,Y,L,BV} proving that any $C^1$-stably weakly shadowable volume-preserving diffeomorphism displays a dominated splitting $E\oplus F$ (Theorem~\ref{T1}). Moreover, both $E$ and $F$ are volume-hyperbolic. With respect to the flow setting the literature is absent on exploring this type of shadowing. We begin were to develop this concept proving similar results (Theorem~\ref{T2}).

\end{subsection}

\begin{subsection}{Basic definitions for the discrete-time case}

Let $M$ be a $d$-dimensional ($d\geq 2$) Riemannian closed and connected manifold and  let $d(\cdot,\cdot)$ denotes the distance on $M$ inherited by the Riemannian structure. We endow $M$ with a volume-form (cf. ~\cite{Mo}) and let $\mu$ denote the Lebesgue measure related to it. Let $\text{Diff}_\mu^{\,\,1}(M)$ denote the set of volume-preserving diffeomorphisms defined on $M$, i.e. those diffeomorphisms such that $\mu(B)=\mu(f(B))$ for any $\mu$-measurable subset $B$. Consider this space endowed with the $C^1$ Whitney topology. The Riemannian inner-product induces a norm $\|\cdot\|$ on the tangent bundle $T_x M$. We will use the usual uniform norm of a bounded linear map $A$ given by $\|A\|=\sup_{\|v\|=1}\|A\cdot v\|$.

Fix some diffeomorphism $f\in \text{Diff}_\mu^{\,\,1}(M)$. Given $\delta>0$, we say that a sequence of points $\{x_i\}_{i\in\mathbb{Z}}\subset M$ is a \emph{$\delta$-pseudo-orbit of $f$} if $d(f(x_i),x_i)<\delta$ for all $i\in\mathbb{Z}$. We say that a sequence of points $\{x_i\}_{i\in\mathbb{Z}}\subset M$ is \emph{weakly $\epsilon$-shadowed by the $f$-orbit of $x$} if $\{x_i\}_{i\in\mathbb{Z}}\subset B({f^i(x),\epsilon})$ for all $i\in\mathbb{Z}$. For an $f$-invariant closed set $\Lambda$ we say that $f|_\Lambda$ has the \emph{weak shadowing property} if for every $\epsilon>0$, there exists $\delta>0$ such that for any $\delta$-pseudo-orbit $\{x_i\}_{i\in\mathbb{Z}}\subset \Lambda$ of $f$, there exists $x\in M$ such that the $f$-orbit of $x$ weakly $\epsilon$-shadows $\{x_i\}_{i\in\mathbb{Z}}$, i.e., $\{x_i\}_{i\in\mathbb{Z}}\subset B({f^i(x),\epsilon})$ for all $i\in\mathbb{Z}$. The diffeomorphism $f$ has the \emph{weak shadowing property} if $M=\Lambda$ in the above definition.

Let $U\subset M$ be a compact set, $f$ a diffeomorphism and $\Lambda_f(U):={\cap}_{n\in\mathbb{Z}}f^n(U)$ the invariant maximal set associated to $U$ and $f$. Take $g$ close to $f$ and let $\Lambda_g(U):={\cap}_{n\in\mathbb{Z}}g^n(U)$ be the continuation of the set $\Lambda_f(U)$. We say that $f|_{\Lambda}$ is \emph{$C^1$-stably weakly shadowing} (in $U$) if there exists a compact neighborhood $U$ of $\Lambda$ such that
\begin{itemize}
\item $\Lambda=\Lambda_f(U)$ and
\item there is a $C^1$-neighborhood $U_f\cap \text{Diff}_\mu^{\,\,1}(M)$ of $f$, such that any $g\in U_f$ has the weakly shadowing property on $\Lambda_g(U)$.
\end{itemize}
Along this paper we deal with the weakly shadowing property on the whole manifold $M$.

Recall, that a set $\Lambda\subseteq M$ is \emph{transitive} if it has a dense orbit. A diffeomorphism $f$ is said to be transitive if $M$ is a transitive set for $f$. We observe that a transitive diffeomorphism has the weakly shadowing property. Thus, $C^1$-stable transitivity, implies the $C^1$-stable shadowing property.

Given an $f$-invariant set $\Lambda\subseteq M$ we say that $\Lambda$ is \emph{uniformly hyperbolic}  if the tangent vector bundle over $\Lambda$ splits into two $Df$-invariant subbundles $T\Lambda=E^u\oplus E^s$ such that $\|Df|_{E^s}\|\leq 1/2$ and $\|Df^{-1}|_{E^u}\|\leq 1/2$. When $\Lambda=M$ we say that $f$ is \emph{Anosov}. Clearly, there are lots of Anosov diffeomorphisms which are not volume-preserving. We say that an $f$-invariant set $\Lambda\subseteq M$ admits an \emph{$\ell$-dominated splitting} if there exists a continuous decomposition of the tangent bundle $T\Lambda$ into $Df$-invariant subbundles $E$ and $F$ such that $$\|Df^{\ell}(x)|_{F}\|.\|(Df^{\ell}(x)|_{E})^{-1}\|\leq 1/2,$$
in this case we say $E \succ_\ell F$ (i.e. $E$ $\ell$-dominates $F$).

Finally, we say that an $f$-invariant set $\Lambda\subseteq M$ is \emph{uniformly partially hyperbolic}, if we have a splitting $E^s\oplus E^c\oplus E^u$ of $T\Lambda$ such that $E^s$ is uniformly contracting, $E^u$ is uniformly expanding, $E^c \succ E^s$ and  $E^u \succ E^c$. When $M$ is partially hyperbolic for $f$ we say that $f$ is a \emph{partially hyperbolic diffeomorphism}.

We say that the $Df$-invariant splitting $E\oplus C\oplus F$ of $T\Lambda$ where $C\succ E$ and $F\succ E$ is \emph{volume partially hyperbolic} if the volume is uniformly contracted on the bundle $E$ and expanded on the bundle $F$, i.e., there exists $\ell\in\mathbb{N}$ such that $|\det Df^\ell(x)_{|E}|<\frac{1}{2}$ and $|\det Df^{-\ell}(x)_{|F}|<\frac{1}{2}$. When $C$ is trivial we say that $E\oplus F$ is \emph{volume hyperbolic}.

It is proved in \cite[Proposition 0.5]{BDP} that, in the volume-preserving context, the existence of a dominated splitting implies volume-hyperbolicity.

A \emph{periodic orbit} for a diffeomorphism $f$ is a point $p\in M$ such that $f^\pi(p)=p$ where $\pi$ is the least positive integer satisfying the equality.  Given a periodic orbit $p$ of period $\pi$ of a diffeomorphism $f$ we say that $p$:
\begin{itemize}
\item is \emph{hyperbolic} if $Df^{\pi}(p)$ has no norm one eigenvalues;
\item has \emph{trivial real spectrum} if $Df^{\pi}(p)$, has only real eigenvalues of norm one (thus equal to $-1$ or $1$) and there exists $0\leq k\leq n$ such that $1$ has multiplicity $k$ and $-1$ has multiplicity $n-k$.
\end{itemize}
Clearly, if $x$ is a $f$-periodic point with period $\pi$ and  with trivial real spectrum, we can split $T_x M=E^+_x\oplus E^-_x$ such that $Df^{\pi}(x)\colon E^+_x \rightarrow E^+_x$ is the identity map and  $Df^{\pi}(x)\colon E^-_x \rightarrow E^-_x$ is the minus identity map.

\end{subsection}

\begin{subsection}{Basic definitions for the continuous-time case} We present the basic set up for the context of flows. In this setting we assume that the dimension of $M$ is greater than two. Given a $C^r$ ($r\geq 1$) vector field $X\colon M \rightarrow TM$ the solution of the equation $x^\prime=X(x)$ gives rise to a $C^r$ flow, $X^t$; by the other side given a $C^r$ flow we can define a $C^{r-1}$ vector field by considering $X(x)=\frac{d X^t(x)}{dt}\vert_{t=0}$. We say that $X$ is \emph{divergence-free} if its divergence is equal to zero. Note that, by Liouville formula, a flow $X^t$ is volume-preserving if and only if the corresponding vector field, $X$, is divergence-free.

Let $\mathfrak{X}_\mu^r(M)$ denote the space of $C^r$ divergence-free vector fields and we consider the usual $C^r$ Whitney topology on this space.

For $\delta>0,$ we say that
$$\{(x_i, t_i):x_i\in M,
t_i\geq1\}_{i\in\Z}$$ is a $(\delta, 1)$-pseudo orbit of $X$ if
$d(X^{t_i}(x_i), x_{i+1})<\delta$ for all $i\in\Z$.

We say that $X$ has the weak shadowing property if for every
$\epsilon>0$ there is $\delta>0$ such that for any $(\delta,
1)$-pseudo orbit $\{(x_i, t_i)\}_{i\in\Z}$, there is a point
$x\in M$ such that $\{x_i\}_{i\in\Z}\subset B(X^t(x),\epsilon)$, for all $t\in\mathbb{R}$.

We say that $X\in\mathfrak{X}_\mu^1(M)$ is \emph{$C^1$-stably weakly shadowable} (on $M$) if any $Y\in\mathfrak{X}_\mu^1(M)$  sufficiently $C^1$-close to $X$ is also weakly shadowable.

Given a vector field $X$ we denote by $Sing(X)$ the set of \emph{singularities} of $X$, i.e. those points $x\in M$ such that $X(x)=\vec{0}$. Let $R:=M\setminus Sing(X)$ be the set of \emph{regular} points.
Given $x\in R$ we consider its normal bundle $N_{x}=X(x)^{\perp}\subset T_{x}M$ and define the associated \emph{linear Poincar\'{e} flow} by $P_{X}^{t}(x):=\Pi_{X^{t}(x)}\circ DX^{t}({x})$ where $\Pi_{X^{t}(x)}\colon T_{X^{t}(x)}M\rightarrow N_{X^{t}(x)}$ is the projection along the direction of $X(X^{t}(x))$. In the same way as we did in the discrete-time case we define can define uniform hyperbolicity, dominated splitting, partial hyperbolicity and volume-hyperbolicity for the linear Poincar\'e flow $P_X^{t}$, in subsets of $R$ and related to subbundles of the normal bundle $N$. We also observe that, when $\Lambda\subseteq M$ is compact, the (partial) hyperbolicity of the tangent map $DX^t$ on $\Lambda$ implies the (partial) hyperbolicity for the linear Poincar\'e flow of $X$ on $\Lambda$ (see \cite[Proposition 1.1]{D}).

Given a closed orbit $p$ of period $\pi$ of a flow $X^t$ we say that $p$:
\begin{itemize}
\item is \emph{hyperbolic} if $P_X^{\pi}(p)$ has no norm one eigenvalues;
\item has \emph{trivial real spectrum} if $P_X^{\pi}(p)$, has only real eigenvalues of norm one (thus equal to $-1$ or $1$) and there exists $0\leq k\leq n-1$ such that $1$ has multiplicity $k$ and $-1$ has multiplicity $n-1-k$.
\end{itemize}
If $x$ is a $X^t$-periodic point with period $\pi$ and with trivial real spectrum, we can split $N_x =N^+_x\oplus N^-_x$ such that $P_X^{\pi}(x)\colon N^+_x \rightarrow N^+_x$ is the identity map and  $P_X^{\pi}(x)\colon N^-_x \rightarrow N^-_x$ is the minus identity map.

\end{subsection}

\begin{subsection}{Statement of the results and some applications}

As we already said, here we begin by developing the generalized versions of the results in \cite{GSW,Y,L,BV}. Our first result is the following.

\begin{maintheorem}\label{T1}
Let $f\in  \text{Diff}_\mu^{\,\,1}(M)$ be a $C^1$-stably weakly shadowing diffeomorphism. Then, $M$ admits a volume-hyperbolic dominated splitting.
\end{maintheorem}

Notice that $f$ is not necessarily $C^1$-robustly transitive and this is exactly the interesting case because otherwise we could use the arguments on robust transitivity developed in \cite[\S 7]{BDP}. We observe that Bonatti and Viana \cite[\S 6.2]{BonV} build an open subset of partially hyperbolic (but not Anosov) transitive diffeomorphisms on $3$-dimensional manifolds. Since their construction can be made conservative we observe that Theorem ~\ref{T1} is optimal for dimension $n\geq 3$.

We observe that a quite complete construction of these type of behavior, in the volume-preserving context, was done in ~\cite{HHTU}. Actually, the authors in ~\cite{HHTU} build \emph{volume-preserving blenders} which are a prototype example of robust transitive dynamics. We also observe that these results were used to prove recent important results in the volume-preserving setting; namely, in ~\cite{C} it is proved that in the complement of Anosov (volume-preserving) diffeomorphisms we have densely robust heterodimensional cycles and, in ~\cite{HHTU2}, it is proved a new criterium for ergodicity among partial hyperbolic diffeomorphisms with central direction with dimension two, that is, stably ergodic diffeomorphisms are $C^1$-dense among volume-preserving partially hyperbolic diffeomorphisms with two-dimensional center bundle.

As we already said, by ~\cite[Proposition 0.5]{BDP}, we obtain that a dominated splitting $TM=E\oplus F$ with $F\succ E$ implies that $E$ is uniformly volume hyperbolic (contracting) and $F$ is uniformly volume hyperbolic (expanding). Thus, Theorem~\ref{T1} implies the result in ~\cite{L}.

With respect to the three-dimensional case Theorem~\ref{T1} implies that, in the presence of $C^1$-stability of the weak shadowing property we get that $TM=E\oplus F$ has a domination $F\succ E$. Since the dimension of the splitting is constant (see ~\cite{BDV}) one of the subbundles $E$ or $F$ are one-dimensional and thus uniformly hyperbolic. In fact, the diffeomorphism is \emph{coarsely} partially hyperbolic (cf. \cite[pp. 122]{H}). In ~\cite[Corollary 1.7]{H} are presented sufficient conditions to have a (proper) partial hyperbolic volume-preserving diffeomorphism, i.e., a splitting into three one-dimensional sub bundles $E^u\oplus E^c\oplus E^s$ with $E^u$ and $E^s$ uniform hyperbolic expanding and contracting, respectively.

Finally, we present the corresponding versions for the flow context.

\begin{maintheorem}\label{T2}
Let $X\in\mathfrak{X}^1_\mu(M)$ be a $C^1$-stably weakly shadowing vector field. Then, $M$ admits a volume-hyperbolic dominated splitting for the linear Poincar\'e flow.
\end{maintheorem}

As an immediate consequence of previous result and \cite[Proposition 1.1]{D} we obtain:

\begin{corollary}
Let $X\in\mathfrak{X}^1_\mu(M)$ be a $C^1$-stably weakly shadowing vector field and $M$ is three-dimensional. Then $X$ is an Anosov flow.
\end{corollary}

We observe that, as we already mention $C^1$-stability of transitivity implies $C^1$-stability of weak shadowing. Thence, Theorem~\ref{T2} gives a different way of proving the main result in ~\cite{BR0}. We notice that Ferreira (see ~\cite{F2}), using the $C^1$-stability of shadowing was able to obtain global hyperbolicity for divergence-free vector fields.

\end{subsection}

\section{Discrete-time case}

\begin{subsection}{Linear conservative cocycles over large periodic systems}\label{LCC}

To prove that any volume-preserving diffeomorphism $f$ which is $C^1$-stably weakly shadowing does not contains trivial real spectrum  periodic orbits we will use the following very useful volume-preserving version of Frank's Lemma proved in \cite[Proposition 7.4]{BDP}.

\begin{lemma}\label{Franks}
Let $f \in \text{Diff}_\mu^{\,\,1}(M)$ and a $C^1$-neighborhood of $f$, $U(f)$, be given. Then there are $\delta_0>0$ and $U_0(f)$ such that for any $g \in U_0(f)$, a finite set $\{ x_1, x_2,...,x_l \}$, a neighborhood U of $\{ x_1, x_2,..., x_l \}$ and volume-preserving maps $L_i: T_{x_i}M \rightarrow T_{g(x_i)}M$ satisfying $||L_i - Dg(x_i)|| < \delta_0$ for all $1 \leq i \leq l$ there are $\epsilon_0>0$ and $\tilde{g} \in U(f)$ such that
\begin{enumerate}
\item [a)] $\tilde{g}(x)=g(x)$ if $x \in M \backslash U$
\item [b)] $\tilde{g} (x)= \varphi_{g(x_i)} \circ L_i \circ \varphi^{-1}_{x_i}(x)$ if $x \in B(x_i,\epsilon_0)$.
\end{enumerate}

\end{lemma}

Assertion b) implies that $\tilde{g}(x)=g(x)$ if $x \in \{ x_1, x_2,....,x_l \}$ and $D\tilde{g}(x_i)$ is conjugated to $L_i$ via the tangent maps of the local volume-preserving charts.

We start this section by recalling some basic definitions introduced in the paper~\cite[\S 2.3]{BGV}. Let $f \in \text{Diff}_\mu^{\,\,1}(M)$ and consider a set $\Sigma \subseteq M$ which is a countable union of periodic orbits of $f$. A \emph{large periods system} (LPS) is a four-tuple $\mathcal{A}=(\Sigma, f, T{\Sigma}, A)$, where $T{\Sigma}$ is the restriction to $\Sigma$ of the tangent bundle over $M$ and $A \colon  \Sigma \rightarrow SL(n,\mathbb{R})$ is a continuous map, where $SL(n,\mathbb{R})$ stands for the special linear $(n^2-1)$-dimensional Lie group of $n\times n$ matrices with real entries. In fact, for $x \in \Sigma$, $A_x$ is a linear map of $T_x M$ and we identify this space with $\mathbb{R}^{n}$. The quintessential  example of a LPS, and associated to the dynamics of the volume-preserving diffeomorphism, is obtained by taking the so-called \emph{dynamical cocycle} given by $A_x=Df(x)$.  Given a LPS $\mathcal{A}=(\Sigma, f, T_{\Sigma}, A)$ the \emph{cocycle identity} associated to it is given by
\begin{equation}\label{lve}
A^{m+n}_x=A^m_{f^n(x)}\cdot A^n_x,
\end{equation}
where $x\in\Sigma$ and $m,n\in\mathbb{N}$.
The LPS $\mathcal{A}=(\Sigma, f, T_{\Sigma}, A)$ is \emph{bounded} if there exists $K>0$ such that $\|A_x\| \leq K$, for all $x \in \Sigma$. The LPS $\mathcal{A}$ is said to be a \emph{large period system} if the number of orbits of $\Sigma$ with period equal to $\tau$ is finite, for any $\tau >0$.  Since the LPS $\mathcal{A}$ evolves in  $SL(n,\mathbb{R})$ we say that it is \emph{conservative}, in fact, $| \det {A}_x| = 1,\,\, \forall x \in \Sigma$.

A  LPS $\mathcal{B}=(\Sigma, f, T{\Sigma}, B)$ is a \emph{conservative perturbation} of a bounded LPS $\mathcal{A}$ if, for every $\epsilon >0$, $\|A_x - B_x\| < \epsilon$,  up to points $x$ belonging to a finite number of orbits, and $\mathcal{B}$ is conservative.  A bounded LPS $\mathcal{A}$ is \emph{strictly without dominated decomposition} if the only invariant subsets of $\Sigma$ that admit a dominated splitting for $A$ are finite sets.

Let us now present a key result about LPS which is the conservative version of ~\cite[Theorem 2.2]{BGV}.

\begin{theorem}\label{2.2}
Let $\mathcal{A}$ be a conservative, large period and bounded LPS. If $\mathcal{A}$ is strictly without dominated decomposition then there exist a conservative perturbation $\mathcal{B}$ of $\mathcal{A}$ and an infinite set $\Sigma^{\prime} \subset \Sigma$ which is  $f$-invariant such that for every $x \in \Sigma^{\prime}$ the linear map $B^{\pi(x)}_x$ as all eigenvalues real and with the same modulus (thus equal to $1$ or to $-1$).
\end{theorem}

As in ~\cite{BGV} we can also obtain the following more user friendly result:

\begin{corollary}\label{cor}
Given any $K>0$ and $\epsilon>0$, there exist $\pi_0,\ell\in \mathbb{N}$ such that for any conservative, large period and $K$-bounded LPS $\mathcal{A}$, over a periodic orbit $x$ with period $\pi(x)>\pi_0$ we have either
\begin{enumerate}
\item [(i)] that $\mathcal{A}$ has an $\ell$-dominated splitting along the orbit of $x$ or else
\item [(ii)] there exists an $\epsilon$-$C^0$-perturbation $\mathcal{B}$ of $\mathcal{A}$, such
that  $B^{\pi(x)}_x$ has all eigenvalues equal to $1$ and $-1$.
\end{enumerate}
\end{corollary}

All the perturbations which are used in the proof of ~\cite[Theorem 2.2]{BGV} are rotations and directional homotheties, i.e., diagonal linear maps for a fixed basis. They are made in the linear cocycle setting which evolves in the general linear group $GL(n,\mathbb{R})$, and clearly can also be done in $SL(n,\mathbb{R})$ with some additional care. Then, Lemma~\ref{Franks} allows to realize them as perturbations of a fixed volume-preserving diffeomorphism. Thence, the proof given by Bonatti, Gourmelon and Vivier can be carried on to our volume-preserving setting without additional obstructions. As a conclusion we obtain the following crucial result.

\begin{lemma}\label{BGV2}
Let $f \in \text{Diff}_\mu^{\,\,1}(M)$ and fix small $\epsilon_0 >0$. There exist $\pi_0,\ell\in \mathbb{N}$ such that for any periodic orbit $x$ with period $\pi(x)>\pi_0$ we have either
\begin{enumerate}
\item [(i)] that $f$ has an $\ell$-dominated splitting along the orbit of $x$ or else
\item [(ii)] for any neighborhood $U$ of $\cup_{n}f^n(x)$, there exists an $\epsilon$-$C^1$-per\-tur\-ba\-tion $g$ of $f$, coinciding with $f$ outside $U$ and on $\cup_{n}f^n(x)$, and such
that  $Dg^{\pi(x)}(x)$ has all eigenvalues equal to $1$ and $-1$.
\end{enumerate}
\end{lemma}

Despite the fact that previous result is very useful we can improve it in order to be used in the sequel.

\begin{remark}
Let $Dg^{\pi(x)}(x)\colon T_x M\rightarrow T_x M$ be the linear map given in Lemma~\ref{BGV2} (ii). We can assume that $x$ has trivial real spectrum. Actually, we can obtain a proof of this by considering ~\cite[Proposition 3.7]{BGV} where we can obtain $h \in \text{Diff}_\mu^{\,\,1}(M)$ arbitrarily $C^1$-close to $g$ and such that $Dh^{\pi(x)}(x)$ has all eigenvalues real of multiplicity 1, and with different modulus. Moreover, the Lyapunov exponents of $Dh^{\pi(x)}(x)$ can be chosen arbitrarily close to those of $Dg^{\pi(x)}(x)$. Finally, another small and volume-preserving perturbation, can be done in order to preserve the simplicity of the spectrum and with all eigenvalues equal to $1$ or $-1$.
\end{remark}

\end{subsection}

\begin{subsection}{Proof of Theorem~\ref{T1}}

The following result is almost the volume-preserving counterpart of ~\cite[Lemma 3.2]{GSW}. Another fundamental results which guarantee the absent of elliptic behavior among $C^1$-stably weakly shadowable maps can be found in \cite[Main Lemma 1]{BV} and in \cite[Proposition 3.5]{L}. In brief terms next lemma says that fixed/periodic local dynamics is an ingredient against the stability of weak shadowing.

\begin{lemma}\label{main}
Fix some $C^1$-weakly shadowable volume-preserving diffeomorphism $f$ and $\delta_0>0$ such that any $g \in \text{Diff}^1_{\mu} (M)$ $\delta_0$-$C^1$-close to $f$ is also weakly shadowable. Let $U_0(f)$ be given by Lemma ~\ref{Franks} with respect to $U(f)$. Then, for any $g \in U_0(f)$, $g$ does not contains periodic points with trivial real spectrum.
\end{lemma}

\begin{proof}
Let  $\dim M=n$ and let us suppose that there is a volume-preserving $g \in U_0(f)$ that have a period orbit $p$ with all eigenvalues equal to 1 and $-1$. Assume that $p$ is such that $g(p)=p$. Then $Dg(p)$ has $k$ eigenvalues equal to 1 and $n-k$ eigenvalues equal to $-1$ and $T_p M=E^{+}_p \oplus E^{-}_p$ where $E^{+}_p$ corresponds to the subspace of the eigenvalue 1 and $E^{-}_p$ corresponds to the subspace of eigenvalue $-1$.

By Lemma ~\ref{Franks} there are $\epsilon_0>0$ and $\tilde{g} \in U(f) \subset \mathcal{WS}(U)$ such that $\tilde{g}(p)= g(p)=p$ and $\tilde{g}(x)= \varphi_{g(p)} \circ Dg(p) \circ \varphi_p^{-1}(x)$ if $x \in B(p,\epsilon_0)$, reducing $\epsilon_0$ if necessary.

The next computations will be yield in $E^{+}_p(\epsilon_0)$ (the other case is similar using that $\tilde{g}^2$ has all eigenvalues equal to 1).

Since $Dg(p)|_{E^{+}_p}= id$, where $id\colon E^{+}_p \rightarrow E^{+}_p$ the identity map on $E^{+}_p$, there is a small arc $\mathcal{I}_p \subset B(p, \epsilon_0) \cap \varphi_p (E^{+}_p(\epsilon_0))$ with center $p$ such that $\tilde{g}(\mathcal{I}_p)=\mathcal{I}_p$ and $\tilde{g}|_{\mathcal{I}_p}$ is  $id$  (i.e., the \emph{identity map}).

Take $\epsilon_1<\epsilon_0$, $v_1 \in E^{+}_p(\epsilon_1)$ with $\|v_1\|=\epsilon_2=\frac{\epsilon_1}{2}$ and set
$$ \mathcal{I}_p \supset \varphi_p( \{ t.v_1:  t \in [-1,1 ]\}) \cap B(p,\epsilon_1).$$

Put $\epsilon= \frac{\epsilon_1}{5}$ and let $0 < \delta < \epsilon$ be the number of the weak shadowing property of $\tilde{g}$.

Now, we are going to construct a $\delta$-pseudo-orbit $\{x_k \}_{k \in \mathbb{Z}}$ of $\tilde{g}$ in $\mathcal{I}_p$ which cannot be weakly $\epsilon$-shadowed by any $\tilde{g}$-true orbit of a point in $M$.

We take a finite sequence $\{ w_k \}_{k=0}^T$ in $E^{+}_p(\epsilon_1)$ for some $T>0$, such that $w_0=O_p$, $w_T=v_1$ and $|w_k - w_{k+1}| < \delta$ for $0 \leq k \leq T-1$. Here the $w_k$ are chosen such that if $w_k=t_k.v_1$, then $|t_k| < |t_{k+1}|$ for $0 \leq k \leq (T-1)$. Finally, define

\begin{itemize}
\item
$x_k=\varphi_p(w_0)$ for $k<0$;
\item
$x_k=\varphi_p(w_k)$ for $0 \leq k \leq T-1$
\item
$x_k=\tilde{g}^{k- T}(\varphi_p(w_T))$ for $k \geq  T$.
\end{itemize}

Then $\{x_k\}_{k \in \mathbb{Z}}$ is a $\delta$-pseudo-orbit  of $\tilde{g}$ in $B(p, \epsilon_2)$ and since $\tilde{g} \in \mathcal{WS}(U)$, there is $y \in M$ weakly $\epsilon$-shadowing $\{ x_k \}_{k \in \mathbb{Z}}.$

The local structure of $\tilde{g}$ in a neighborhood of $\mathcal{I}_p$ in $M$ is the direct product of the  identity map, $\tilde{g}|_{\varphi_p(E^{+}_p(\epsilon_2)) \cap B(p,\epsilon_2)}$ by the minus identity map $\tilde{g}|_{\varphi_p(E^{-}_p(\epsilon_2)) \cap B(p,\epsilon_2)}$.

We may assume that $y \in B(x_0,\epsilon)$.

If $y \in \mathcal{I}_p$ then since  $\tilde{g}^i(y)=y$ for $i \in \mathbb{Z}$, $d(\tilde{g}^i (y), x_T)> \epsilon$ by the choice of $\epsilon$.

If $y \notin \mathcal{I}_p$  then, $\tilde{g}(y)=y$ or else $\tilde{g}^2(y)=y$, for all $i \in \mathbb{Z}$, and $d(\tilde{g}^i (y), x_T)> \epsilon$ by the choice of $\epsilon$. This proves the lemma.
\end{proof}

\begin{proof}(of Theorem~\ref{T1})
Take the Pugh-Robinson residual (general density theorem, see ~\cite{PR}) intersected with Bonatti-Crovisier residual ~\cite{BC} and call it $R$. As $\text{Diff}_\mu^{\,\,1}(M)$ endowed with the $C^1$-topology is a Baire space we can take a sequence of $f_n \in R$ with $f_n\rightarrow f$ (in the $C^1$-topology). Since, $M$ is transitive for each $f_n$, there exists periodic points $p_n$ (for $f_n$ and with period $\pi_n$) such that $\lim \sup_n \cup_i f^i(p_n)=M$ (in the Hausdorff metric sense\footnote{We recall that the Hausdorff distance between two compact subsets $A,B\subseteq M$ is given by $d_H(A,B)=\max\{\sup_{y\in B}d(y,A),\sup_{x\in A}d(x,B)\}$.}). Clearly, $\pi_n\rightarrow +\infty$. Define
\begin{equation}\label{sigma0}
\Sigma=\bigcup_{n\in\mathbb{N}}\{p_n,f_n(p_n),...,f_n^{\pi_n-1}(p_n)\}.
\end{equation}
Notice that $\overline{\Sigma}=M$.

We now define a $\mathcal{A}=(\Sigma,g,T\Sigma,A)$; $\Sigma$ is defined in (\ref{sigma0}), $T\Sigma=T_x M$ where $x\in\Sigma$, $g(f^i_n(p_n))=f^{i+1}_n(p_n)$ for any $i\in\{0,1,...,\pi_n-1\}$ and $A(g^i_n(p_n))=Df_n(f^i_n(p_n))$.

Of course that we are in the presence of a LPS. By Corollary~\ref{cor} there exists a uniform dominated splitting over $\Sigma$. Since $\overline{\Sigma}=M$ and the dominated splitting extends to the closure (cf. \cite{BDV}) we obtain that $M$ has a (uniform) dominated splitting with respect to the LPS  $\mathcal{A}=(\Sigma,g,T\Sigma,A)$.

Realizing dynamically the cocycle,  by Theorem~\ref{2.2} and Lemma~\ref{main}, we obtain that $f$, $C^1$-stably weakly shadowable, has a dominated splitting over $M$. Finally, ~\cite[Proposition 0.5]{BDP} guarantees the volume-hyperbolic statement  and the theorem is proved.

\end{proof}

\end{subsection}

\section{Continuous-time case}

\begin{subsection}{Linear traceless differential systems over large periodic systems}

We begin by recalling some definitions introduced in~\cite{BGV}, in \S~\ref{LCC} and first developed for flows in ~\cite{BR0}. Let $X \in \mathfrak{X}_\mu^1(M)$ and consider a set $\Sigma \subset M$ which is a countable union of closed orbits of $X^t$. A \emph{Linear Differential System} (LDS) is a four-tuple $\mathpzc{A}=(\Sigma, X^t, N_{\Sigma}, A)$, where $N_{\Sigma}$ is the restriction to $\Sigma$ of the normal bundle of $X$ over $M\setminus Sing(X)$ and $A \colon  \Sigma \rightarrow GL(n-1,\mathbb{R})$ is a continuous map. In fact, for $x \in \Sigma$, $A_x$ is a linear map of $N_x$ and we identify this space with $\mathbb{R}^{n-1}$. The natural LDS associated to the dynamics of the vector field is obtained by taking $A_x=\Pi \circ DX_x$.

Given $\mathpzc{A}=(\Sigma, X^t, N_{\Sigma}, A)$ the \emph{linear variational equation} (or \emph{equation of first variations}) associated to it is given by:
\begin{equation}\label{lve}
\dot{u}(t,x)=A(X^t(x)) \cdot u(t,x).
\end{equation}

The matriciant (or solution) of the system (\ref{lve}) with initial condition $u(0,x)=id$ is, for each $t$ and $x$, a linear map $\Phi_{A}^t(x) \colon N_x \rightarrow N_{X^t(x)}$. We call the map $A$ the \emph{infinitesimal generator} of $\Phi_ {A}$ . It is easy to see that  $\Phi_{A}^t(x)=P_X^t(x)$ when the infinitesimal generator is $\Pi \circ DX$. We say that $\mathpzc{A}=(\Sigma, X^t, N_{\Sigma}, A)$ is \emph{bounded}, if there exists $K>0$ such that $\|A_x\| \leq K$, for all $x \in \Sigma$. The LDS $\mathpzc{A}$ is said to be a \emph{large period system} if the number of orbits of $\Sigma$ with period less or equal to $\tau$ is finite, for any $\tau >0$. We say that the LDS $\mathpzc{A}$ is \emph{traceless}\footnote{We observe that, if for a given $X\in\mathfrak{X}^1_\mu(M)$, we have $Sing(X)=\emptyset$, then there exists $\rho\colon M\rightarrow \mathbb{R}$ with $\rho(x)=\|X(x)\|^{-1}$ such that $\rho X(x)$, in (\ref{MAP}), satisfies $| \det\Phi_{A}^t(x) |=1$. Now, by (\ref{LF}) we get $\int_0^t tr(A(X^s(x)))ds=0$. } (or \emph{conservative}) if, for all $x \in \Sigma$, we have
\begin{equation}\label{MAP}
| \det\Phi_{A}^t(x) |\|X(X^t(x))\| = \|X(x)\|.
\end{equation}
It follows from Liouville's formula that
\begin{equation}\label{LF}
\det\Phi_{A}^t(x) = e^{\int_0^t tr(A(X^s(x)))ds}.
\end{equation}

A  LDS $\mathpzc{B}=(\Sigma, X^t, N_{\Sigma}, B)$  is a \emph{traceless perturbation} of a bounded LDS $\mathpzc{A}$ if, for every $\epsilon >0$, $\|A_x - B_x\| < \epsilon$,  up to points $x$ belonging to a finite number of orbits, and $\mathpzc{B}$ is conservative. From (\ref{LF}) it follows that $\mathpzc{B}$ is traceless if and only if $tr(B)=tr(A)$.

Gronwall's inequality gives that $$\|\Phi_{A}^t(x)-\Phi_{B}^t(x)\| \leq \exp(K|t|)\|A_x - B_x\|.$$
In particular $\Phi_{B}^1$ is a perturbation of $\Phi_{A}^1$ in the sense introduced in~\cite{BGV} and in \S~\ref{LCC} for the discrete case. A bounded LDS $\mathpzc{A}$ is \emph{strictly without dominated decomposition} if the only invariant subsets of $\Sigma$ that admit a dominated splitting for $\Phi_{A}^t$ are finite sets.

Now we present a result about LDS which is the flow version of Theorem ~\ref{2.2}.

\begin{theorem}\label{3.2}
Let $\mathpzc{A}$ be a traceless, large period and bounded LDS. If $\mathpzc{A}$ is strictly without dominated decomposition then there exist a traceless perturbation $\mathpzc{B}$ of $\mathpzc{A}$ and an infinite set $\Sigma^{\prime} \subset \Sigma$ which is  $X^t$-invariant such that for every $x \in \Sigma^{\prime}$ the linear map $\Phi_{B}^{\pi(x)}(x)$ as all eigenvalues real and with the same modulus (thus equal to $1$ or to $-1$).
\end{theorem}

Like in Theorem~\ref{2.2}, once again the perturbations are made in the linear traceless differential systems and Franks' Lemma for volume-preserving vector fields (Lemma~\ref{mpl}) allows to realize them as perturbations of a fixed volume-preserving flow.

Let us introduce now a useful concept. Fix $X\in\mathfrak{X}_\mu^1(M)$ and $\tau>0$. A \emph{one-parameter area-preserving linear family} $\{A_t\}_{t\in \mathbb{R}}$ associated to $\{X^t(p);\,\, t \in [0, \tau]\}$ is defined as follows:
\begin{itemize}
\item $A_t\colon N_p \rightarrow N_p$ is a linear map, for all $t\in \mathbb{R}$,
\item $A_t=id$, for all $t\leq 0$, and $A_t=A_{\tau}$, for all $t\geq \tau$,
\item $A_t\in SL(n, \mathbb{R})$, and
\item the family $A_t$ is $C^\infty$ on the parameter $t$.
\end{itemize}

The following result, proved in ~\cite[Lemma 3.2]{BR0} is now stated for $X \in \mathfrak{X}_\mu^1(M)$ instead of $X \in \mathfrak{X}_\mu^4(M)$ because of the improved \emph{smooth $C^1$ pasting lemma} proved in ~\cite[Lemma 5.2]{BR2}.

\begin{lemma}\label{mpl}
Given $\epsilon>0$ and a vector field $X \in \mathfrak{X}_\mu^1(M)$ there exists $\xi_0=\xi_0(\epsilon,X)$ such that $\forall \tau \in [1,2]$, for any periodic point $p$ of period greater than $2$, for any sufficient small flowbox $\mathcal{T}$ of $\{X^t(p);\,\, t \in [0, \tau]\}$ and for any one-parameter linear family $\{A_t\}_{t \in [0, \tau]}$ such that $\|  A_t^{\prime} A_t^{-1}\|<\xi_0$, $\forall t \in [0, \tau]$, there exists $Y \in \mathfrak{X}_\mu^1(M)$ satisfying the following properties
\begin{enumerate}
\item $Y$ is $\epsilon$-$C^1$-close to $X$;
\item $Y^t(p)=X^t(p)$, for all $t \in \mathbb{R}$;
\item $P_Y^\tau(p)=P_X^\tau(p) \circ A_{\tau}$, and
\item $Y|_{\mathcal{T}^c}\equiv X|_{\mathcal{T}^c}$.
\end{enumerate}
\end{lemma}

In overall, we obtain the following result:

\begin{lemma}\label{BGV3}
Let $X\in \mathfrak{X}^1_{\mu} (M)$ and fix a small $\epsilon_0 >0$. There exist $\pi_0,\ell\in \mathbb{N}$ such that for any closed orbit $x$ with period $\pi(x)>\pi_0$ we have either
\begin{enumerate}
\item [(i)] that $P_X^t$ has an $\ell$-dominated splitting along the orbit of $x$ or else
\item [(ii)] for any neighborhood $U$ of $\cup_{t}X^t(x)$, there exists an $\epsilon$-$C^1$-per\-tur\-ba\-tion $Y$ of $X$, coinciding with $X$ outside $U$ and on $\cup_{t}X^t(x)$, and such that  $P_Y^{\pi(x)}(x)$ has all eigenvalues equal to $1$ and $-1$.
\end{enumerate}
\end{lemma}

\end{subsection}

\begin{subsection}{Set us free of singularities}
In order to rule out singularities in the context of $C^1$-stably weak shadowable volume-preserving flows we will recall some useful results.

The first one was proved in (\cite[Lemma 3.3]{BR0}).

\begin{lemma}\label{linear}
Let $\sigma$ be a singularity of $X \in \mathfrak{X}_\mu^1(M)$. For any $\epsilon >0$ there exists  $Y \in  \mathfrak{X}_\mu^{\infty}(M)$, such that $Y$ is $\epsilon$-$C^1$-close to $X$ and  $\sigma$ is a linear hyperbolic singularity of $Y$.
\end{lemma}

The second one, was proved in ~\cite[Proposition 4.1]{V} generalizing Doering theorem in \cite{D}. Observe that, in our context, the singularities of hyperbolic type are all saddles.

\begin{proposition}\label{Vivier}
If $Y \in \mathfrak{X}^1(M)$ admits a linear hyperbolic singularity of saddle-type, then the linear Poincar\'e flow of $Y$ does not admit any dominated splitting  over $M\setminus Sing(Y)$.
\end{proposition}

Finally, since by Poincar\'e recurrence, any $X \in \mathfrak{X}_\mu^1(M)$ is chain transitive, the following result is a direct consequence of \cite{B}.

\begin{proposition}\label{Mixing}
In $\mathfrak{X}_\mu^1(M)$ chain transitive flows equal topologically mixing flows in a $C^1$-residual subset.
\end{proposition}

The following theorem is proved borrowing some arguments in \cite[Theorem 15]{AR}.

\begin{theorem}\label{sing}
If $X\in\mathfrak{X}_\mu^1(M)$ is $C^1$-stably weak shadowable vector field, then $X$ has no singularities.
\end{theorem}
\begin{proof}
Let $X\in\mathfrak{X}_\mu^1(M)$ be a $C^1$-stably weak shadowable vector field and fix a small $C^1$ neighborhood $\mathcal{U}\subset \mathfrak{X}_\mu^1(M)$ of $X$. The proof is by contradiction. Assume that $Sing(X)\not= \emptyset$. Using Lemma~\ref{linear}, there exists $Y\in \mathcal{U}$ with a linear saddle-type singularity $\sigma\in Sing(Y)$. By Proposition~\ref{Mixing}, there exist $Z_n\in \mathfrak{X}_\mu^1(M)$ $C^1$-close to $Y$ which is topologically mixing. We can find $W_n\in \mathfrak{X}_\mu^1(M)$ $C^1$-close to $Z_n$ having a $W_n$-closed orbit $p_n$ such that the Hausdorff distance between $M$ and $\cup_t W_n^t(p_n)$ is less than $1/n$.

Now we consider jointly Lemma~\ref{main2} and Lemma~\ref{BGV3} obtaining that $P_{W_n}^t$ is $\ell$-dominated over the $W_n$-orbit of $p_n$ where $\ell$ is uniform on $n$. Since $W_n$ converges in the $C^1$-sense to $Y$ and $\lim \sup_n \cup_t W_n^t(p_n)=M$ we obtain that $M\setminus Sing(Y)$ has an $\ell$-dominated splitting which contradicts Proposition~\ref{Vivier}.
\end{proof}

\end{subsection}

\begin{subsection}{Proof of Theorem~\ref{T2}}

The following result is the flow counterpart of Lemma~\ref{main}.

\begin{lemma}\label{main2}
Fix some $C^1$-weakly shadowable volume-preserving vector field $X\in\mathfrak{X}^1_{\mu} (M)$. Then, any $Y \in \mathfrak{X}^1_{\mu} (M)$ sufficiently $C^1$-close to $X$ does not contains closed orbits with trivial real spectrum.
\end{lemma}

\begin{proof}

Let $X \in \mathfrak{X}_\mu^1(M)$ be a volume preserving vector field and $U(X)$ a $C^1$-neighbourhood of X in $\mathfrak{X}_\mu^1(M)$ where the weak shadowing property holds. Let $p$ be a closed orbit of $X$ with period $\pi>1$ and $U_p$ a small neighbourhood of $p$ in $M$. Let us also assume that all eigenvalues of $P^{\pi}_X(p)$, the linear Poincar\'e flow in $p$ are $-1$ and $1$.

Now, we transpose our objects to the euclidian space using the volume-preserving charts given by Moser's theorem (\cite{Mo}). Then, there exists a smooth conservative change of coordinates $\varphi_p\colon T_pM \rightarrow U_p$ such that $\varphi_p(\vec{0})=p$. Let $f_X: \varphi_p(N_p) \rightarrow \Sigma$ stands for the Poincar\'e map associated to $X^t$, where $\Sigma$ denotes the Poincar\'e section through $p$ and take $\mathcal{V}$ a $C^1$-neighbourhood of $f_X$.

By the careful construction of a \emph{local linearized} divergence-free vector field done in ~\cite{BR1} we take $\mathcal{T}$ a small flowbox of $\{X^t(p)\colon t \in [0,\tau] \}$, $\tau < \pi$ where we have that there are a (local linear) divergence-free vector field $Z \in U(X)$, $f_Z \in \mathcal{V}$ and a small $\epsilon_0>0$ such that

$$ f_Z(x)= \left\{
      \begin{array}{ll}
        \varphi_p \circ P_X^{\pi} \circ \varphi_p^{-1}(x), & \hbox{$x \in B_{\epsilon_0}(p) \cap \varphi_p(N_p)$;} \\
        f_X(x), & \hbox{$x \notin B_{4\epsilon_0}(p) \cap \varphi_p(N_p)$.}
      \end{array}
    \right. $$

The next computations will be yield in $N^{+}_p(\epsilon_0)$. Take $v \in N^{+}_p(\epsilon_1)$, $\epsilon_1 <\epsilon_0$ with $||v||=\epsilon_2 = \frac{\epsilon_1}{2}$ and set $\mathcal{I}_p=\{s.v\colon 0 \leq s \leq 1 \}$.

Fix $0< \epsilon < \frac{\epsilon_2}{2}$ and let $0 < \delta < \epsilon$ be the number of the weak shadowing property of $Z^t$. Now, we are going to construct a $(\delta, 1)$-pseudo-orbit of $Z^t$ belonging to $\varphi_p(\mathcal{I}_p)$ which cannot be weakly $\epsilon$-shadowed by any true orbit  $y \in M$.

We take a finite sequence $\{w_k \}_{k=0}^T \in N^{+}_p(\epsilon_1)$ for some $T>0$, such that $w_0=0_p$ and $w_T=v$ and $|w_k-w_{k+1}|< \delta$ for $0 \leq k \leq T-1$. Here $w_k$ are chosen such that if $w_k=s_k.v$ then $s_k<s_{k+1}$ for $0 \leq k \leq T-1$. Finally, we define,

\begin{itemize}
\item
$x_k=\varphi_p(w_0)$, $t_k=\pi$ for $k<0$;
\item
$x_k=\varphi_p(w_k)$, $t_k=\pi$ for\footnote{Observe that we are considering that the return time at the transversal section is the same and equal to $\pi$. Clearly, it is not exactly equal to $\pi$, however it is as close to $\pi$ as we want by just squeezing the flowbox.} $0 \leq k \leq T-1$
\item
$x_k=f_Z^{k-T}(\varphi_p(w_T))$,  $t_k=\pi$ for $k \geq  T$.
\end{itemize}

Then $\{(x_i,t_i)\}_{i \in \mathbb{Z}}$ is a $(\delta,1)$-pseudo orbit of $Z^t$ in $B(p, \epsilon_2)$ and since $Z$ is weakly shadowable, there is $y \in M$ such that $\{x_i\}_{i \in \mathbb{Z}} \subset B(Z^t(y), \epsilon)$, $\forall t \in \mathbb{R}.$

We may assume that $y \in B(x_0, \epsilon) \cap \varphi_p(N_{p,\epsilon})$. Note that, in $N^{+}_p(\epsilon_1)$ we have

\begin{eqnarray*}
d(x_0, x_T)&=&d(\varphi_p(w_0),\varphi_p(w_T))=d(\varphi_p(\vec0),\varphi_p(v))\\
&=&d(p,\varphi_p(v))\approx \|v\|=\epsilon_2>2\epsilon,
\end{eqnarray*}
where $\approx$ says that $d(p,\varphi_p(v))$ is very close to $d(\vec0,v)=\|v\|$, because $\varphi_p$ is very close to be an isometry when $\epsilon$ is close to zero.

On the other hand, since $Z$ is weakly shadowable, we have that for some $\iota= n \pi$ we have,

$$d(x_0,x_T) \leq d(x_0,y)+d(y,x_T)=d(x_0,y)+d(Z^{\iota}(y),x_T)<2\epsilon,$$
which is a contradiction and the lemma follows.
\end{proof}

\begin{proof}(of Theorem~\ref{T2})
The proof goes as the one did in Theorem~\ref{T1}. Take the Pugh-Robinson residual (general density theorem, see ~\cite{PR}) intersected with the residual in ~\cite{B} and call it $R$. As $\mathfrak{X}^1_{\mu} (M)$ endowed with the $C^1$-topology is a Baire space we can take a sequence of $X_n \in R$ with $X_n\rightarrow X$ (in the $C^1$-topology). Since $M$ is transitive, for each $X_n$, there exists closed orbits $p_n$ (for $X_n$ and with period $\pi_n$) such that $\lim \sup_n \cup_t X^t(p_n)=M$ in the Hausdorff metric sense. Clearly, $\pi_n\rightarrow +\infty$. Define
\begin{equation}\label{sigma}
\Sigma=\bigcup_{n\in\mathbb{N}}\{X_n^t(p_n)\colon 0\leq t \leq \pi_n\}.
\end{equation}
Notice that $\overline{\Sigma}=M$.

We now define a LDS $\mathpzc{A}=(\Sigma,Y^t,N_\Sigma,A)$; $\Sigma$ is defined in (\ref{sigma}), $N_\Sigma=N_x$ where $x\in\Sigma$. We define a one-parameter map by $Y^t(X^r_n(p_n))=X^{t+r}_n(p_n)$ for any $r\in[0,\pi_n]$. Observe that $Y^t$ is a flow; clearly 
\begin{itemize}
\item $Y^0(X^r_n(p_n))=X^{r}_n(p_n)$ and
\item $Y^{t+s}(X^r_n(p_n))=X^{t+s+r}_n(p_n) = Y^t(X^{s+r}_n(p_n))=Y^t(Y^s(X^{r}_n(p_n)))$.
\end{itemize}

Finally, we define the linear action on $N_\Sigma$ by 
$$\Phi_A^t(Y^r_n(p_n))=P_{X_n}^t(X^r_n(p_n)).$$ 
Since $X_n$ are divergence-free, the map $A$ is traceless.

This time we are in the presence of a LDS. Reasoning analogously with the discrete-time case, using Theorem~\ref{3.2}, we obtain that there exists a uniform dominated splitting over $\Sigma$. Since $\overline{\Sigma}=M$ and $Sing(Y)=\emptyset$, the dominated splitting extends to the closure and we obtain that $M$ has a (uniform) dominated splitting with respect to the LDS $\mathcal{A}=(\Sigma,Y^t, N_\Sigma,A)$.

Finally, we realize dynamically the traceless LDS,  by Theorem~\ref{3.2} and Lemma~\ref{mpl}, we obtain that $X^t$, $C^1$-stably weakly shadowable, has a dominated splitting over $M$ for $P_X^t$. The volume-hyperbolicity can be obtained as in ~\cite[Proposition 0.5]{BDP}.

\end{proof}

\end{subsection}

\section{Conclusion and future directions}

In overall, as it was expected, shadowing is a property that allows us to go further than weak shadowing. However, surprisingly, does not take us as far as we would expect when compared to weak shadowing. In fact, in low dimensions, the $C^1$-stability of the shadowing property implies global hyperbolicity exactly as the $C^1$-stability of the weak shadowing property. Moreover, even in the multidimensional case, both properties share several properties, because both avoid the presence of singularities (for flows), both display a dominated splitting in the whole manifold and both have a hyperbolic behavior; the shadowing have hyperbolicity and the weak shadowing volume-hyperbolicity.

We believe that the techniques developed in the present paper should be useful to prove that volume-preserving dynamical systems which exhibits $C^1$-shadowing-like properties display a weak form of hyperbolicity. Actually, it is our guess that the results in \cite{AR} can be generalized not only for our class but also for dissipative flows in higher dimensions. Thus, systems with $C^1$-stability of the  average shadowing property and also the asymptotic average shadowing property should also have a dominated splitting. Furthermore, the $C^1$-stability of weak shadowing for dissipative flows is not studied yet, and our results should enlighten the solution for that problem.

\section*{Acknowledgements}

MB and SV were partially supported by National Funds through FCT - ``Funda\c{c}\~{a}o para a Ci\^{e}ncia e a Tecnologia", project \\PEst-OE/MAT/UI0212/2011.
 ML was supported by Basic Science Research Program through the National
Research Foundation of Korea(NRF) funded by the Ministry of
Education, Science and Technology, Korea (No. 2011-0007649).

\end{document}